\documentclass[11pt]{article}

\usepackage{times,amssymb,amsmath,exscale,array,latexsym}
\usepackage{graphicx}
\usepackage{epsfig}
\usepackage{color}
\definecolor{marin}{rgb}   {0.,   0.3,   0.7}
\definecolor{rouge}{rgb}   {0.8,   0.,   0.}
\definecolor{sepia}{rgb}   {0.8,   0.5,   0.}
\usepackage[colorlinks,citecolor=marin,linkcolor=rouge,
            bookmarksopen,
            bookmarksnumbered
           ]{hyperref}

\newcommand{\e}{\ensuremath{\mathrm{e}}}

\newcommand{\psis}{{\mathcal S}}

\numberwithin{equation}{section}

\newcommand{\QED}{\mbox{}\hfill \raisebox{-0.2pt}{\rule{5.6pt}{6pt}\rule{0pt}{0pt}}
          \medskip\par}

\begin{document}
\title{High order integrators obtained by linear combinations of symmetric-conjugate compositions
}

\author{F. Casas$^{1}$, A. Escorihuela-Tom\`as$^{2}$ \\[2ex]
$^{1}$ {\small\it Departament de Matem\`atiques and IMAC, Universitat Jaume I, 12071-Castell\'on, Spain}\\{
\small\it Email: Fernando.Casas@mat.uji.es}\\[1ex]
$^{2}$ {\small\it Departament de Matem\`atiques, Universitat Jaume I, 12071-Castell\'on, Spain}\\{
\small\it Email: alescori@uji.es}\\[1ex]
}


%
\maketitle

\begin{abstract}

A new family of methods  involving complex coefficients for the numerical integration of differential equations is presented and analyzed.
They are constructed as linear combinations of symmetric-conjugate compositions obtained from a basic time-symmetric integrator of order
$2n$ ($n \ge 1$). The new integrators are of order $2(n+k)$, $k=1,2,\ldots$, and preserve time-symmetry up to order $4n+3$ when
applied to differential equations with real vector fields. 
If in addition the system is Hamiltonian and the basic scheme is symplectic, then they also preserve symplecticity
up to order $4n+3$. We show that these integrators are well suited for a parallel implementation, thus improving their efficiency.
Methods up to order 10 based on a 4th-order integrator are built and tested in comparison with other standard procedures to increase the order of a basic scheme.

\vspace*{1cm}


\end{abstract}\bigskip

\noindent AMS numbers: 65L05, 65P10, 37M15


\section{Introduction}
\label{sec.1}
Composition methods constitute a standard tool to construct high-order numerical integrators for the initial value problem
\begin{equation}  \label{ivp}
  \dot{x} = f(x), \qquad x(t_0) = x_0 \in \mathbb{R}^d,
\end{equation}
in particular when the vector field $f$ possesses some qualitative property whose preservation by numerical approximations is deemed relevant
\cite{blanes16aci,hairer06gni}. Let 
$\mathcal{S}_h^{[2n]}$ denote a $2n$-th order method, so that $\mathcal{S}_h^{[2n]}(x_0) = \varphi_h(x_0) + \mathcal{O}(h^{2n+1})$,  
where $x(h) = \varphi_h(x_0)$ is the exact solution of Eq. (\ref{ivp}) for a time step $h$. Then, if the coefficients $\alpha_1, \alpha_2, \ldots, \alpha_s$
satisfy some algebraic conditions, the composition of the basic scheme with step sizes $\alpha_1 h, \alpha_2 h, \ldots, \alpha_s h$, i.e.,
\begin{equation}  \label{compo.2b}
  \psi_h= \mathcal{S}_{\alpha_1 h}^{[2n]} \circ \mathcal{S}_{\alpha_{2} h}^{[2n]} \circ \cdots \circ 
  \mathcal{S}_{\alpha_{s-1} h}^{[2n]} \circ \mathcal{S}_{\alpha_s h}^{[2n]}
\end{equation}
is a new method of higher order $2n + m$ \cite{blanes08sac}. 
If in particular $f$ is Hamiltonian and $\mathcal{S}_h^{[2n]}$ is symplectic, then the composition method
(\ref{compo.2b}) is also symplectic \cite{hairer06gni}. In general, any geometric property the basic method has in common with the exact solution is 
still shared by the higher-order scheme (\ref{compo.2b}) if this property is preserved by composition \cite{mclachlan02sm}. 
Moreover, suppose $\mathcal{S}_h^{[2n]}$ is time-symmetric,
namely, it satisfies
\[
  \mathcal{S}_h^{[2n]} \circ \mathcal{S}_{-h}^{[2n]} = \mathrm{id},
\]  
where $\mathrm{id}$ is the identity map, for any $h$. Then, method (\ref{compo.2b}) is also time-symmetric if the composition is left-right palindromic, i.e., $\alpha_{s+1-j} = \alpha_{j}$,
$j=1,2,\ldots$.

A well known class of composition methods is obtained by applying the triple-jump procedure \cite{suzuki90fdo,yoshida90coh}: 
\begin{equation} \label{eq:tj}
  \mathcal{S}_h^{[2n+2]} = \mathcal{S}_{\alpha_1 h}^{[2n]} \circ \mathcal{S}_{\alpha_2 h}^{[2n]} \circ \mathcal{S}_{\alpha_1 h}^{[2n]},
\end{equation}
with
\begin{equation} \label{real.1}
  \alpha_1 = \frac{1}{2 - 2^{1/(2n+1)}}, \qquad \alpha_2 = 1 - 2 \alpha_1,
\end{equation}
is a new method of order $2n+2$.
The same technique can be applied again to $\mathcal{S}_h^{[2n+2]}$, so that one can construct recursively time-symmetric methods of any order $2n+2k$, 
$k=1,2,\ldots$.

When constructing high-order composition methods, real coefficients $\alpha_1, \ldots, \alpha_s$ are not the only option, however. In fact, the unavoidable 
existence of negative $\alpha_j$ in (\ref{compo.2b}) when the order is higher than two \cite{blanes05otn,goldman96nno,sheng89slp,suzuki91gto} 
typically imposes 
stability restrictions on the step size. This occurs in particular when
Eq. (\ref{ivp}) is the outcome of a parabolic differential equation discretized in space. In that case, considering complex coefficients with positive real part 
is also a valid alternative \cite{castella09smw,hansen09hos}. Even for problems where the presence of some $\alpha_j < 0$ is not particularly troublesome, 
composition methods with
complex coefficients have also been proposed and analyzed from the preservation of properties viewpoint \cite{blanes10smw,blanes21osc,chambers03siw}.

In the particular case of the triple-jump composition (\ref{eq:tj}), in addition to the real solution (\ref{real.1}), the complex one with the smallest phase is
\begin{equation} \label{cco1}
  \alpha_1 = \frac{\e^{i \pi / (2n+1)}}{2^{1/(2n+1)} - 2 \e^{i \pi / (2n+1)}}, \qquad \alpha_2 = 1 - 2 \alpha_1,
\end{equation}
and the resulting method has  in fact smaller truncation errors than its real counterpart (\ref{real.1}). 
If the basic scheme is time-symmetric and of order 2, then time-symmetric
methods up to order 14 with coefficients having positive real part are possible by applying this technique \cite{blanes13oho}.

The order can be raised by one instead with the simplest composition \cite{bandrauk91ies,suzuki91gto}
\begin{equation} \label{sc.1}
  \psi_h^{[2n+1]} = \mathcal{S}_{\alpha_1 h}^{[2n]} \circ \mathcal{S}_{\alpha_{2} h}^{[2n]}
\end{equation}
if
\[
  \alpha_1 = \bar{\alpha}_2 = \frac{1}{2} + \frac{i}{2} \, \frac{\sin \frac{2 \ell +1}{2n+1} \pi}{1 + \cos \frac{2 \ell +1}{2n+1} \pi} \qquad \mbox{ for }
  \qquad -n \le \ell \le n-1.
\]  
The choice $\ell = 0$ gives the solution with the smallest phase, which we denote by $\gamma^{[2n]}$:
\begin{equation}  \label{5c}
   \alpha_1 = \gamma^{[2n]} := \frac{1}{2} + \frac{i}{2} \, \frac{\sin \frac{\pi}{2n+1}}{1 + \cos \frac{\pi}{2n+1}}, \qquad n = 1,2, \ldots
\end{equation}
When the vector field $f$ in (\ref{ivp}) is real, then $x_1 = \psi_h^{[2n+1]}(x_0)$ is complex, and so it is quite natural to project $x_1$ on the real axis and
proceed to the next step only with $\Re(x_1)$. This is equivalent of course to integrating with the scheme
\begin{equation}  \label{R.1}
   R_{h}^{(1)} =\frac{1}{2} \left( \psi_{h}^{[2n+1]} + \overline{\psi}_{h}^{[2n+1]} \right). 
\end{equation}
Method (\ref{R.1}) is not time-symmetric, even when $\mathcal{S}_h^{[2n]}$ is. Nevertheless, it has been shown in \cite{casas21cop} that $R_{h}^{(1)}$ is
\emph{pseudo-symmetric} of order $4n+3$, in the sense that
\[
 R_h^{(1)} \circ R_{-h}^{(1)} = \mathrm{id} + \mathcal{O}(h^{4n+4})
\]  
if the vector field $f$ in (\ref{ivp}) is real. If in addition $f$ is Hamiltonian and $\mathcal{S}_h^{[2n]}$ is symplectic, then $R_{h}^{(1)}$ is
also   \emph{pseudo-symplectic} of order $4n+3$. In other words, projecting $\psi_h^{[2n+1]}$ at each integration step leads to a numerical method that
preserves geometric properties of the exact solution up to an order that is much higher than the order of the method itself. Pseudo-symplectic
integrators have been previously considered in the literature, both in the context of Runge--Kutta \cite{aubry98psr} and polynomial extrapolation methods
\cite{blanes99eos,chan00eos}.

Moreover, as shown in  \cite{casas21cop},
$R_{h}^{(1)}$ can be taken as the basis of the recursion
\begin{equation}   \label{foit}
  R_{h}^{(k)}=    \frac{1}{2} \left( R_{\gamma^{[2k]} h}^{(k-1)} \circ R_{\bar{\gamma}^{[2k]} h}^{(k-1)}  +  
     R_{\bar{\gamma}^{[2k]} h}^{(k-1)}  \circ R_{\gamma^{[2k]} h}^{(k-1)}  \right), \qquad\qquad  k = 2, 3 \ldots,
\end{equation}
producing methods of order $2(n+k)$, also pseudo-symmetric of order $4n+3$. Here the coefficients $\gamma^{[2k]}$ are given by
Eq. (\ref{5c}). For future reference, we call (\ref{foit}) \textbf{$R$-methods}.

Scheme (\ref{sc.1}) is a particular example of a \emph{symmetric-conjugate} composition. These are composition methods of the form
\begin{equation}  \label{compo.2}
  \psi_h= \mathcal{S}_{\alpha_1 h}^{[2n]} \circ \mathcal{S}_{\alpha_{2} h}^{[2n]} \circ \cdots \circ 
  \mathcal{S}_{\bar{\alpha}_2 h}^{[2n]} \circ \mathcal{S}_{\bar{\alpha}_1 h}^{[2n]},
\end{equation}
i.e., compositions (\ref{compo.2b}) with $\alpha_j \in \mathbb{C}$ and
\[
   \bar{\alpha}_{s+1-j} = \alpha_j, \qquad j=1,2,\ldots.
\]   
Methods of this class, 
as shown in \cite{blanes21osc}, possess remarkable preservation properties when considering its real part, 
\[
    \Re(\psi_h) = \frac{1}{2} \left( \psi_{h} + \overline{\psi}_{h} \right).
\]
 In particular, if one takes a time-symmetric 2nd-order scheme as the basic method and the coefficients $\alpha_1, \alpha_2, \ldots$
are chosen in such a way that $\psi_h$ is of order $2n-1$,  then  $\Re(\psi_h)$ is of order $2n$ and pseudo-symmetric of order $4n-1$ 
  when the vector field $f$ in  (\ref{ivp}) is real.
  If in addition $f$ is a (real) Hamiltonian vector field and $\psis_{h}^{[2]}$ is a symplectic integrator, then $ \Re(\psi_h)$ is pseudo-symplectic of order $4n-1$.

Since taking the real part of a symmetric-conjugate method is just a very special linear combination, it is quite natural to ask what happens when one
considers a more general linear combination of symmetric-conjugate compositions and their complex-conjugate, $\psi_h^{(j)}$, $\overline{\psi}_h^{(j)}$: is it possible
to construct new methods of higher order whereas still preserving time-symmetry (and symplecticity) up to the order prescribed by the composition 
$\psi_h^{(j)}$? If yes, how the new methods are built? Addressing these questions is precisely the subject of the present paper. In doing so, we present a new family
of schemes of increasingly higher order well adapted for implementation in a parallel environment, 
requiring less computational effort than the $R$-methods (\ref{foit}) but with the same qualitative properties.

If we denote for simplicity the symmetric-conjugate composition (\ref{compo.2}) by its sequence of coefficients,
\[
    \psi_h^{(j)}= (\alpha_1, \alpha_{2}, \ldots, \alpha_{s-1}, \alpha_s),
\]
with $\bar{\alpha}_{s+1-j} = \alpha_j$, these new schemes have the basic structure
\begin{equation} \label{T.methods}
  T_{h}^{(k)} = \frac{1}{2^k} \sum_{j=1}^{2^{k-1}} \Big( (\alpha_{j_{2^k}},  \ldots,  \alpha_{j_1})  + \mbox{c.c.} \Big)
\end{equation}
and are of order $2(n+k) \le 4n+3$ and pseudo-symmetric of order $4n+3$. We designate them as \textbf{$T$-methods}.

\section{Construction of the family of $T$-methods}
\label{sec.2}

In this section we construct the new family of integrators $T_{h}^{(k)}$ and show explicitly that they are of order $2(n+k)$ and pseudo-symmetric of order $4n+3$ for
$k=1,2,3$. The same procedure can be formally extended to any $k > 3$. The analysis is based on the Lie formalism applied to the series of differential
operators associated to the integrators.

\subsection{Series of differential operators}

As is well known, given a time-symmetric integrator $\mathcal{S}_{h}^{[2n]}$ of order $2n  \ge 2$
one can associate a series of linear operators $\exp(Y(h))$ so that
\[
    g(\mathcal{S}_{h}^{[2n]}(x)) = \exp(Y(h))[g](x) 
\]
for all functions $g$ \cite{blanes08sac}, with 
\[
     Y(h) = h Y_1 + h^{2n+1} Y_{2n+1} + h^{2n+3} Y_{2n+3} + \cdots.
\]     
Here $Y_k$ are certain operators depending on the particular method and, for consistency, $Y_1 = F$, where $F$ is
the Lie derivative associated with $f$:
\begin{equation}  \label{eq.1.1b}
  F = \sum_{i\ge1} \, f_i(x) \, \frac{\partial }{\partial x_i}.
\end{equation}
The composition (\ref{compo.2b}) then has the associated series
\begin{equation} \label{eq.2.4}
  \Psi(h) =   \exp(Y(h \alpha_s)) \, \exp(Y(h \alpha_{s-1})) \, \cdots \,
  \exp(Y(h \alpha_{2})) \, \exp(Y(h \alpha_1)),
\end{equation}
which can be  formally written as $\Psi(h) = \exp(V(h))$ by repeated application of the
Baker--Campbell--Hausdorff formula, with
\[
  V(h) = h F +
h^{2n+1} V_{2n+1} + h^{2n+2} V_{2n+2} + \cdots.
\]
Here $V_{2n+1}, V_{2n+2}, \ldots$ are linear combinations of Lie brackets involving the operators
$Y_1, Y_{2n+1}, Y_{2n+3},\ldots$ \cite{mclachlan02sm}. 
In the particular case of a symmetric-conjugate composition (\ref{compo.2}),
terms $V_{2k}$ in $V(h)$ of even powers in $h$ are pure imaginary, whereas terms $V_{2k+1}$ are real \cite{blanes21osc}.

 For a consistent symmetric-conjugate composition (\ref{compo.2}), i.e., verifying
\begin{equation} \label{consistency}
  \sum_{j=1}^s \alpha_j =  \alpha_1 + \alpha_2 + \cdots + \bar{\alpha}_2 +  \bar{\alpha}_1  = 1,
\end{equation}  
 we get explicitly
\begin{equation} \label{vh}
  V(h) =h  E_{1,1} + h^{2n} \sum_{j \ge 0} h^{2j+1} 
       \sum_{k=1}^{\ell_{2j+1}} \mu_{2n+2j+1,k} E_{2j+1,k} + i \, h^{2n} \sum_{j \ge 1} h^{2j} \sum_{k=1}^{\ell_{2j}} \sigma_{2n+2j,k} E_{2j,k},
\end{equation}
where $\mu_{n,k}$,  $\sigma_{n,k}$ are homogeneous real polynomials of degree $n$ in
the coefficients $\alpha_l$, $l=1, \ldots, s$, and $E_{n,k}$ are elements $Y_j$ and independent Lie brackets involving these operators.
In particular
\[
  \mu_{2n+2j+1,1}=\sum_{l=1}^s \alpha_l^{2(n+j)+1}, \quad\quad j\geq 0
\] 
and
\[
E_{1,1} = Y_1, \quad E_{2n+2n+2j+1,1} = Y_{2n+2j+1}, \quad 
E_{2n+2j,1} = [E_{1,1}, E_{2n+(2j-1),1}], \quad  j=1,2, \ldots
\]

\subsection{Linear combinations of symmetric-conjugate compositions}

Let us now consider the linear combination
\begin{equation} \label{phi}
  \phi_h = \frac{1}{2k} \sum_{j=1}^k \left( \psi_h^{(j)} + \overline{\psi}_h^{(j)} \right),
\end{equation}
where each $\psi_h^{(j)}$ is a consistent symmetric-conjugate composition of the form (\ref{compo.2}) with
different coefficients $\alpha_k^{(j)}$. Then, clearly, $\phi_h$ has 
\begin{equation} \label{def.Phi}
  \Phi(h) \equiv \frac{1}{2k} \sum_{j=1}^k \left( \Psi^{(j)}(h) + \overline{\Psi}^{(j)}(h) \right) = 
      \frac{1}{2k} \sum_{j=1}^k  \left( \e^{V_j(h)} + \e^{\overline{V}_j(h)} \right)
\end{equation}      
as the associated series of operators, where each $V_j(h)$ is of the form (\ref{vh}). Now, by following the same approach as in \cite{casas21cop}, we 
express $\Phi(h)$ as
\[
  \Phi(h) =  \frac{1}{2k}  \e^{\frac{h}{2} F}  \, \sum_{j=1}^k  \left( \e^{W_j(h)} + \e^{\overline{W}_j(h)} \right) \, \e^{\frac{h}{2} F},
\]
where
\begin{equation} \label{expreW}
\begin{aligned}
  & W_j(h)  =  h^{2n+1} \mu_{2n+1,1}^{(j)} E_{2n+1,1} + i \, h^{2n+2} \sigma_{2n+2,1}^{(j)} E_{2n+2,1} \\
    & \qquad + h^{2n+3} \left( \mu_{2n+3,1}^{(j)} E_{2n+3,1}  + \Big( \mu_{2n+3,2}^{(j)} + \frac{1}{24} \mu_{2n+1,1}^{(j)} \Big) \, E_{2n+3,2} \right) \\
    & \qquad + i \, h^{2n+4} \left( \sigma_{2n+4,1}^{(j)} E_{2n+4,1} + \Big( \sigma_{2n+4,2}^{(j)} + \frac{1}{24} \sigma_{2n+2,1}^{(j)} \Big) \, E_{2n+4,2} \right)  \\
   & \qquad + h^{2n+5} \left( \mu_{2n+5,1}^{(j)} E_{2n+5,1} + \Big( \mu_{2n+5,2}^{(j)} + \frac{1}{24} \mu_{2n+3,1}^{(j)} \Big) \, E_{2n+5,2} \, + \right. \\
      & \qquad\qquad\quad \left. \Big(  \mu_{2n+5,3}^{(j)} +  \frac{1}{24} \mu_{2n+3,2}^{(j)} + \frac{1}{1920} \mu_{2n+1,1}^{(j)} \Big) \, E_{2n+5,3} \right)  \\   
     & \qquad + i \, h^{2n+6} \left( \sigma_{2n+6,1}^{(j)} E_{2n+6,1} + \Big( \sigma_{2n+6,2}^{(j)} + \frac{1}{24} \sigma_{2n+4,1}^{(j)} \Big) \, 
            E_{2n+6,2} + \right.  \\
    & \qquad\qquad\quad \left. \Big(  \sigma_{2n+6,3}^{(j)} +  \frac{1}{24} \sigma_{2n+4,2}^{(j)} + \frac{1}{1920} \sigma_{2n+2,1}^{(j)} \Big) \, E_{2n+6,3} \right) \\
    & \qquad + \mathcal{O}(h^{2n+7}).
 \end{aligned}
 \end{equation}
 Here 
 \[
 \begin{aligned}
  &  E_{2n+3,2} = [E_{1,1}, E_{2n+2,1}], \quad  E_{2n+4,2} = [E_{1,1}, E_{2n+3,2}], \quad  E_{2n+5,2} = [E_{1,1}, E_{2n+4,1}], \\
  &  E_{2n+6,2} = [E_{1,1}, E_{2n+5,1}], \quad  E_{2n+5,3} = [E_{1,1}, E_{2n+4,2}], \quad  E_{2n+6,2} = [E_{1,1}, E_{2n+5,2}].
\end{aligned}
\] 
 This is done by applying the symmetric Baker--Campbell--Hausdorff formula to each product $\e^{-\frac{h}{2} F} \, \e^{V_j(h)} \, \e^{-\frac{h}{2} F}$.
From (\ref{expreW}), a straightforward calculation shows that
\[
\begin{aligned}
 & (W_j + \overline{W}_j)^2 = 4 h^{4n+2}  (\mu_{2n+1,1}^{(j)})^2 \, E_{2n+1,1}^2 + \mathcal{O}(h^{4n+4}) \\
 & W_j^2 + \overline{W}_j^2 = 2  h^{4n+2}  (\mu_{2n+1,1}^{(j)})^2 \,  E_{2n+1,1}^2 + \mathcal{O}(h^{4n+4}).
\end{aligned}
\]
Therefore, 
\[
  \frac{1}{2} \left( \e^{W_j} + \e^{\overline{W}_j} \right) - \e^{\frac{1}{2} (W_j + \overline{W}_j)} = \frac{1}{4} (W_j^2 + \overline{W}_j^2) -
  \frac{1}{8} (W_j + \overline{W}_j)^2 + \cdots = \mathcal{O}(h^{4n+4})
\]
and $\Phi(h)$ can also be written as
 \[
  \Phi(h) =  \frac{1}{k}   \, \sum_{j=1}^k \e^{\frac{h}{2} F}  \, \e^{\frac{1}{2} (W_j(h) + \overline{W}_j(h))}   \,  \e^{\frac{h}{2} F} 
  + \mathcal{O}(h^{4n+4}).
\]
In consequence, each term in $\phi_h$ is time-symmetric up to terms $h^{4n+3}$, with independence of the polynomials $\mu_{k,l}^{(j)}$,
$\sigma_{k,l}^{(j)}$, since the sum $W_j(h) + \overline{W}_j(h)$ only contains odd powers of $h$.

On the other hand, one has
\[
\begin{aligned}
  &  \frac{1}{2k} \sum_{j=1}^k \left( \e^{W_j} + \e^{\overline{W}_j} \right) - \exp \left(\frac{1}{2k} \sum_{j=1}^k (W_j + \overline{W}_j) \right) \\
  & =  \frac{1}{4k} \sum_{j=1}^k (W_j^2 + \overline{W}_j^2) - \frac{1}{8k^2} \left( \sum_{j=1}^k (W_j + \overline{W}_j) \right)^2 + \cdots \\
  & =  h^{4n+2} \,  \frac{1}{2k} \left( \sum_{j=1}^k (\mu_{2n+1}^{(j)})^2 - \frac{1}{k} \Big( \sum_{j=1}^k \mu_{2n+1,1}^{(j)} \Big)^2 \right) \, E_{2n+1,1}^2 +
   \mathcal{O}(h^{4n+4}), 
\end{aligned}   
\]
so that it is also true that
\begin{equation} \label{psi5}
  \Phi(h) =   \exp\left(\frac{h}{2} F\right)  \exp \left( \frac{1}{2k} \sum_{j=1}^k  \big(W_j(h) + \overline{W}_j(h) \big) \right) \exp \left(\frac{h}{2} F \right) +
  \mathcal{O}(h^{4n+2}).
\end{equation}

\subsection{Order conditions}

It is thus possible to obtain the order conditions for the method $\phi_h$ in (\ref{phi}) by analyzing just the exponent of the central term in (\ref{psi5}).
From (\ref{expreW}) it follows that
\[
\begin{aligned}
 & \frac{1}{2} \sum_{j=1}^k  \left(W_j(h) + \overline{W}_j(h) \right) = h^{2n+1}  \, c_{2n+1,1} \, E_{2n+1,1} + 
   h^{2n+3} c_{2n+3,1} E_{2n+3,1} \\
   & \qquad + h^{2n+3}  \left(c_{2n+3,2} + \frac{1}{24}  c_{2n+1,1} \right) \, E_{2n+3,2} +h^{2n+5} \, c_{2n+5,1} \, E_{2n+5,1} \\ 
  & \qquad    + h^{2n+5}  \left(c_{2n+5,2} + \frac{1}{24} c_{2n+3,1} \right) \,  E_{2n+5,2} \\
 & \qquad +   h^{2n+5} \left( c_{2n+5,3} + \frac{1}{24}  c_{2n+3,2} + \frac{1}{1920} c_{2n+1,1} \right) \, E_{2n+5,3} + \mathcal{O}(h^{2n+7}),
 \end{aligned}
 \]
with
\[
\begin{aligned}
 & c_{2n+1,1} =  \sum_{j=1}^k  \mu_{2n+1,1}^{(j)}, \qquad c_{2n+3,1} = \sum_{j=1}^k  \mu_{2n+3,1}^{(j)}, \qquad c_{2n+3,2}= \sum_{j=1}^k  \mu_{2n+3,2}^{(j)} \\
 & c_{2n+5,1} = \sum_{j=1}^k  \mu_{2n+5,1}^{(j)}, \qquad c_{2n+5,2 }= \sum_{j=1}^k  \mu_{2n+5,2}^{(j)}, \qquad  c_{2n+5,3} =\sum_{j=1}^k  \mu_{2n+5,3}^{(j)} 
\end{aligned}
\] 
 
In consequence, for consistent compositions $\psi_h^{(j)}$, $j=1,\ldots, k$, the conditions to be satisfied so that $\phi_h$ is a method of
order $r$ are the following:
\begin{itemize} 
 \item $r=2n+2$: $\, c_{2n+1,1} = 0$
 \item $r=2n+4$: $\, c_{2n+1,1} = c_{2n+3,1} = c_{2n+3,2} = 0$
 \item $r=2n+6$: $\, c_{2n+1,1} = c_{2n+3,1} = c_{2n+3,2} = c_{2n+5,1} = c_{2n+5,2} = c_{2n+5,3} = 0 $
\end{itemize}

\subsection{New schemes}

Once identified the relevant order conditions, our
next goal is to solve these equations with the minimum number of basic schemes in the compositions $\psi_h^{(j)}$ and the minimum value of $k$
in the linear combination (\ref{phi}).

\paragraph{Order $r=2n+2$.}

One needs to solve two equations to get a method $\phi_h$ of order $2n+2$: consistency and $c_{2n+1,1}=0$. These can be satisfied by
taking $k=1$ and the simplest composition $\psi_h =  \mathcal{S}_{\alpha_1 h }^{[2n]} \circ \mathcal{S}_{\bar{\alpha}_1 h}^{[2n]}$, in which case one has
\[
   \alpha_1 + \bar{\alpha}_1 = 1, \qquad \alpha_1^{2n+1} + \bar{\alpha}_1^{2n+1} = 0.
\]
In other words, we recover the composition (\ref{sc.1}) and the $R$-method (\ref{R.1}). Our first $T$-method (\ref{T.methods}) is thus
\begin{equation} \label{eq.t1}
  T_{h}^{(1)}   =  \frac{1}{2} \Big( (\gamma^{[2n]}, \bar{\gamma}^{[2n]}) + ( \bar{\gamma}^{[2n]}, \gamma^{[2n]}) \Big) 
\end{equation}
or in more detail
\[
  T_{h}^{(1)}   =  \frac{1}{2} \Big( \psis_{\gamma^{[2n]} h}^{[2n]} \circ  \psis_{\bar{\gamma}^{[2n]} h}^{[2n]} + 
     \psis_{\bar{\gamma}^{[2n]} h}^{[2n]}  \circ \psis_{\gamma^{[2n]} h}^{[2n]}  \Big).
\]      

\paragraph{Order $r=2n+4$.}

Now we have to solve 3 order conditions in addition to consistency for the compositions $\psi_h^{(j)}$ involved. As before, one 
could take in principle $k=1$. In that case, the minimum number of basic maps in $\psi_h^{(1)}$ is 4, just to have enough parameters 
to satisfy the order conditions. It turns out, however, that there are no solutions with the required symmetry
$\alpha_4 = \bar{\alpha}_1$, $\alpha_3 = \bar{\alpha}_2$. In fact, if we take
\[
   \psi_h^{(1)} = (\bar{\alpha}_1, \bar{\alpha}_2, \alpha_2, \alpha_1), \qquad \mbox{ with } \qquad
     \alpha_1 = \bar{\gamma}^{[2n+4]} \bar{\gamma}^{[2n+2]}, \quad   \alpha_2 = \bar{\gamma}^{[2n+4]} \gamma^{[2n]}, 
\]
then  $\mu_{2n+1,1}^{(1)} = \mu_{2n+3,1}^{(1)} = 0$, but  $\mu_{2n+3,2}^{(1)} \ne 0$. On the other hand, if we take
\[
 \psi_h^{(2)} = (\bar{\alpha}_2, \bar{\alpha}_1, \alpha_1, \alpha_2)
\]
with the \emph{same} values of $\alpha_1$, $\alpha_2$ as before, then $\mu_{2n+3,2}^{(2)} = - \mu_{2n+3,2}^{(1)}$, whereas
still verifying that $\mu_{2n+1,1}^{(2)} = \mu_{2n+3,1}^{(2)} = 0$. In consequence, by combining both compositions,
\[
   \phi_h = \frac{1}{4} \left( \psi_h^{(1)} + \overline{\psi}_h^{(1)} + \psi_h^{(2)} + \overline{\psi}_h^{(2)} \right),
\]
one gets a method of order $2n+4$ and pseudo-symmetric of order $4n+3$. This corresponds to our second $T$-method, which
reads explicitly
\begin{eqnarray} \label{eq.t2}
  T_{h}^{(2)}  & = & \frac{1}{4} \Big( (\gamma^{[2n+2]} \gamma^{[2n]}, \gamma^{[2n+2]} \bar{\gamma}^{[2n]}, \bar{\gamma}^{[2n+2]} \gamma^{[2n]},
         \bar{\gamma}^{[2n+2]} \bar{\gamma}^{[2n]} ) \nonumber \\
   &    & \quad  +   (\gamma^{[2n+2]} \bar{\gamma}^{[2n]}, \gamma^{[2n+2]} \gamma^{[2n]}, \bar{\gamma}^{[2n+2]} \bar{\gamma}^{[2n]},
         \bar{\gamma}^{[2n+2]} \gamma^{[2n]} ) \\
   &   & \quad  +   (\bar{\gamma}^{[2n+2]} \bar{\gamma}^{[2n]}, \bar{\gamma}^{[2n+2]} \gamma^{[2n]}, \gamma^{[2n+2]} \bar{\gamma}^{[2n]},
         \gamma^{[2n+2]} \gamma^{[2n]} )   \nonumber \\  
   &    & \quad +  (\bar{\gamma}^{[2n+2]} \gamma^{[2n]}, \bar{\gamma}^{[2n+2]} \bar{\gamma}^{[2n]}, \gamma^{[2n+2]} \gamma^{[2n]},
         \gamma^{[2n+2]} \bar{\gamma}^{[2n]} ) \Big).   \nonumber
\end{eqnarray}
Again, the coefficients $\gamma^{[2m]}$ are given by Eq. (\ref{5c}).

\paragraph{Order $r=2n+6$.}
A total of 7 equations (including consistency) have to be solved in this case, so that we take a symmetric-conjugate composition
involving $s=8$ basic maps,
\[
   \psi_h^{(1)} = (\alpha_1, \alpha_2, \alpha_3, \alpha_4, \bar{\alpha}_4, \bar{\alpha}_3, \bar{\alpha}_2, \bar{\alpha}_1).
\]
With the choice
\[
\begin{aligned}
 &  \alpha_1 =  \gamma^{[2n+4]} \gamma^{[2n+2]} \gamma^{[2n]}, \qquad\quad \alpha_2 =  \gamma^{[2n+4]} \gamma^{[2n+2]} \bar{\gamma}^{[2n]}, \\
 & \alpha_3 =  \gamma^{[2n+4]} \bar{\gamma}^{[2n+2]} \gamma^{[2n]}, \qquad\quad \alpha_4 = \gamma^{[2n+4]} \bar{\gamma}^{[2n+2]} \bar{\gamma}^{[2n]}
\end{aligned} 
\]
it turns out that conditions $c_{2n+1}=c_{2n+3,1}=c_{2n+5,1}=0$ are automatically satisfied. By following the same approach as before, we permute
the position of the coefficients and take the composition
\[
 \psi_h^{(2)} = (\alpha_2, \alpha_1, \alpha_4, \alpha_3, \bar{\alpha}_3, \bar{\alpha}_4, \bar{\alpha}_1, \bar{\alpha}_2).
\]
Then, one has $\mu_{2n+3,2}^{(2)} = - \mu_{2n+3,2}^{(1)}$, so that $\psi_h^{(1)} + \psi_h^{(2)}$ leads to a method of order $2n+4$. More
composition have to be incorporated, however, in order to verify conditions $c_{2n+5,2} =0$ and $c_{2n+5,3} =0$. The former 
is accomplished  by both sums  $\psi_h^{(1)} + \psi_h^{(4)}$ and $\psi_h^{(2)} + \psi_h^{(3)}$, where
\[
\begin{aligned}
  &   \psi_h^{(3)} = (\alpha_3, \alpha_4, \alpha_1, \alpha_2, \bar{\alpha}_2, \bar{\alpha}_1, \bar{\alpha}_4, \bar{\alpha}_3) \\
  &   \psi_h^{(4)} = (\alpha_4, \alpha_3, \alpha_2, \alpha_1, \bar{\alpha}_1, \bar{\alpha}_2, \bar{\alpha}_3, \bar{\alpha}_4), 
\end{aligned}
\]
but the later is satisfied only by adding up the four compositions. In summary, the linear combination
\[
   \frac{1}{4} \Re(  \psi_h^{(1)} +  \psi_h^{(2)} +  \psi_h^{(3)} +  \psi_h^{(4)} ) 
\]   
leads to a method  of order $2n+6$, denoted as $T_h^{(3)}$. More explicitly,
\begin{equation}  \label{deft3}
       T_h^{(3)} = \frac{1}{8} \Big(  \psi_h^{(1)} +  \psi_h^{(2)} +  \psi_h^{(3)} +  \psi_h^{(4)}  + \overline{\psi}_h^{(1)} + \overline{\psi}_h^{(2)}
       + \overline{\psi}_h^{(3)} + \overline{\psi}_h^{(4)} \Big).
\end{equation}

\

The same procedure can be carried out in general, although more order conditions (and consequently more compositions involving more
basic maps) have to be dealt with. This class of methods can be represented in a convenient way as follows.
If we introduce the matrix of coefficients
\[
   \Gamma_{2n} := \frac{1}{2} \left(  \begin{array}{cc}
   					\gamma^{[2n]}  &  \bar{\gamma}^{[2n]}  \\
					\bar{\gamma}^{[2n]}  &  \gamma^{[2n]} 
				\end{array}  \right)
\]
then, according with the previous results, method $T_{h}^{(1)}$ (of order $2n+2$) can be represented by $\Gamma_{2n}$, 
\[
  T_{h}^{(1)} \leadsto   \Gamma_{2n}, 
\]   
whereas
$T_{h}^{(2)}$ (of order $2n+4$) can be associated with the matrix 
\[
     \Gamma_{2n+2} \otimes \Gamma_{2n} = \frac{1}{4} \left(  \begin{array}{cccc}
     	\gamma^{[2n+2]} \gamma^{[2n]}  &  \gamma^{[2n+2]} \bar{\gamma}^{[2n]}  &  \bar{\gamma}^{[2n+2]} \gamma^{[2n]}  &  \bar{\gamma}^{[2n+2]} \bar{\gamma}^{[2n]}  \\
     	\gamma^{[2n+2]} \bar{\gamma}^{[2n]}  &  \gamma^{[2n+2]} \gamma^{[2n]}  &  \bar{\gamma}^{[2n+2]} \bar{\gamma}^{[2n]}  & \bar{\gamma}^{[2n+2]} \gamma^{[2n]}     \\
     	\bar{\gamma}^{[2n+2]} \gamma^{[2n]}  &  \bar{\gamma}^{[2n+2]} \bar{\gamma}^{[2n]}  &  \gamma^{[2n+2]} \gamma^{[2n]}  & \gamma^{[2n+2]} \bar{\gamma}^{[2n]}     \\
     	\bar{\gamma}^{[2n+2]} \bar{\gamma}^{[2n]} & \bar{\gamma}^{[2n+2]} \gamma^{[2n]}  &   \gamma^{[2n+2]} \bar{\gamma}^{[2n]}  &  \gamma^{[2n+2]} \gamma^{[2n]}  
	   \end{array}  \right),
\]
in the sense that each file of $\Gamma_{2n+2}$ corresponds to a particular symmetric-conjugate composition entering into the formulation of $T_{h}^{(2)}$. 
We can write analogously
\[
   T_{h}^{(2)} \leadsto \Gamma_{2n+2} \otimes \Gamma_{2n},
\]   
and moreover
\[
   T_{h}^{(3)} \leadsto \Gamma_{2n+4} \otimes (\Gamma_{2n+2} \otimes \Gamma_{2n}).
\]   
In general, the coefficients in  the $T$-method of order $r=2n+2k$ are distributed according with the pattern
\[
   T_{h}^{(k)} \leadsto \Gamma_{2(n+ k-1)} \otimes (  \Gamma_{2(n+ k -2)}  \otimes \cdots \otimes (\Gamma_{2n+2} \otimes \Gamma_{2n}) \cdots ). 
\]

\section{Numerical examples}

We illustrate next the behavior of some of the previously constructed 
$T$-methods on a pair of numerical examples.  The first one (the 2-dimensional Kepler problem) allows one
to check preservation properties, whereas the second (a simple diffusion equation) is used as a test of their relative performance.
In all cases we take as basic scheme $\mathcal{S}_h^{[2n]}$ the 4th-order ($n=2$) time-symmetric splitting method
\begin{equation}  \label{s4sim}
  \psis_{h}^{[4]} = \varphi_{b_1 h}^{[b]} \circ \,  \varphi_{a_1 h}^{[a]} \circ \, \varphi_{b_2 h}^{[b]} \circ \,  \varphi_{a_2 h}^{[a]} \circ \, 
  \varphi_{b_3 h}^{[b]} \circ \,  \varphi_{a_2 h}^{[a]} \circ \, \varphi_{b_2 h}^{[b]} \circ \,  \varphi_{a_1 h}^{[a]} \circ \, \varphi_{b_1 h}^{[b]} 
\end{equation}
with coefficients
 \begin{eqnarray} \label{P4S9bis} 
b_1 &=& 0.060078275263542357774 - 0.060314841253378523039 \, i, \\  
a_1  &=& 0.18596881959910913140, \nonumber \\  
b_2  &=& 0.27021183913361078161 + 0.15290393229116195895 \, i,  \nonumber \\  
a_2  &=& 0.31403118040089086860, \nonumber \\  
 b_3 &=& 0.33941977120569372122 - 0.18517818207556687181 \, i,  \nonumber 
 \end{eqnarray}
previously considered in \cite{blanes13oho}. This integrator is intended for Eq. (\ref{ivp}) when $f$ can be decomposed as $f(x) = f_a(x) + f_b(x)$ in such a 
way that each sub-problem
\[
   \dot{x} = f_a(x), \qquad\qquad  \dot{x} = f_b(x), 
\]
with $x(0) = x_0$, has solution $x(t)= \varphi_{t}^{[a]}(x_0)$, and $x(t)= \varphi_{t}^{[b]}(x_0)$, respectively.   

\

The implementation of all the integrators has been done in Python 3.7 running on Debian GNU/Linux 10 and the operations with complex arithmetics 
have been coded using the complex class of the numpy library.

\paragraph{Kepler problem.}

The Hamiltonian function for the planar two-body problem reads
\begin{equation}   \label{eq.HamKepler}
   H(q,p) = T(p) + V(q) = 
	\frac{1}{2} p^T p - \mu \frac{1}{r}.
\end{equation}
Here $q=(q_1,q_2)$, $p=(p_1,p_2)$, $r= \|q\|$, $\mu=GM$, $G$ is the gravitational constant and 
$M$ is the sum of the masses of the two bodies. The
corresponding equations of motion are then
\[
   \dot{q}_i = \frac{\partial H}{\partial p_i}  = p_i, \qquad\qquad   \dot{p}_i = -\frac{\partial H}{\partial q_i} = -\mu \frac{q_i}{r^3}, \qquad i=1,2.
\]
Taking $\mu=1$ and initial conditions
\begin{equation}\label{eq.1.12}
  q_1(0) = 1- e, \quad q_2(0) = 0, \quad p_1(0) = 0, \quad p_2(0) = \sqrt{\frac{1+e}{1-e}},
\end{equation}
the resulting trajectory is an ellipse of eccentricity $0 \le e < 1$. In this case $\varphi_{h}^{[a]}$ (respectively, $\varphi_{h}^{[b]}$)
corresponds to the exact solution obtained by integrating the kinetic energy $T(p)$ (resp., potential energy $V(q)$) in 
(\ref{eq.HamKepler}). 

\begin{figure}[!h]
  \begin{center}
      \includegraphics[scale=0.55]{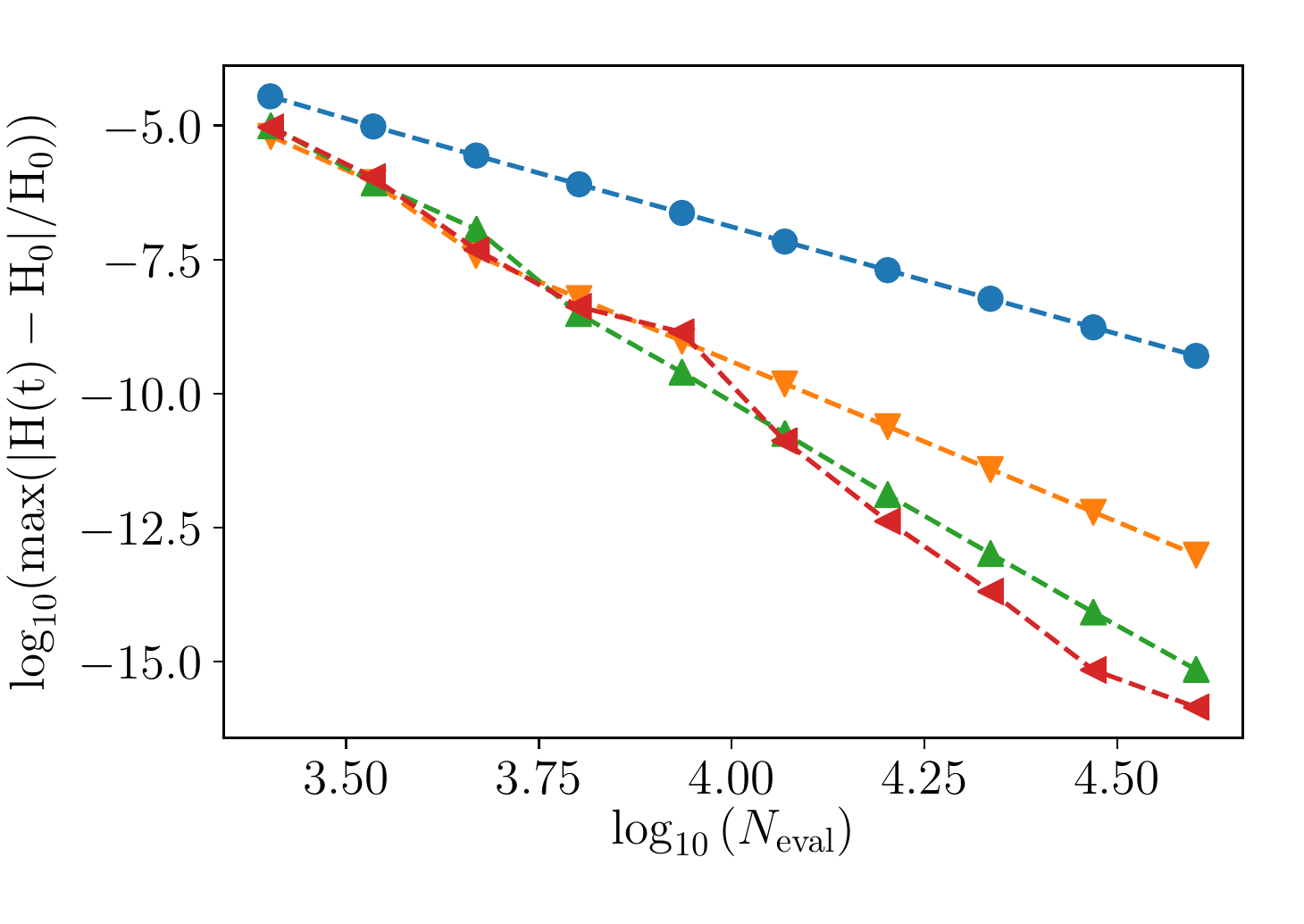}
      \includegraphics[scale=0.55]{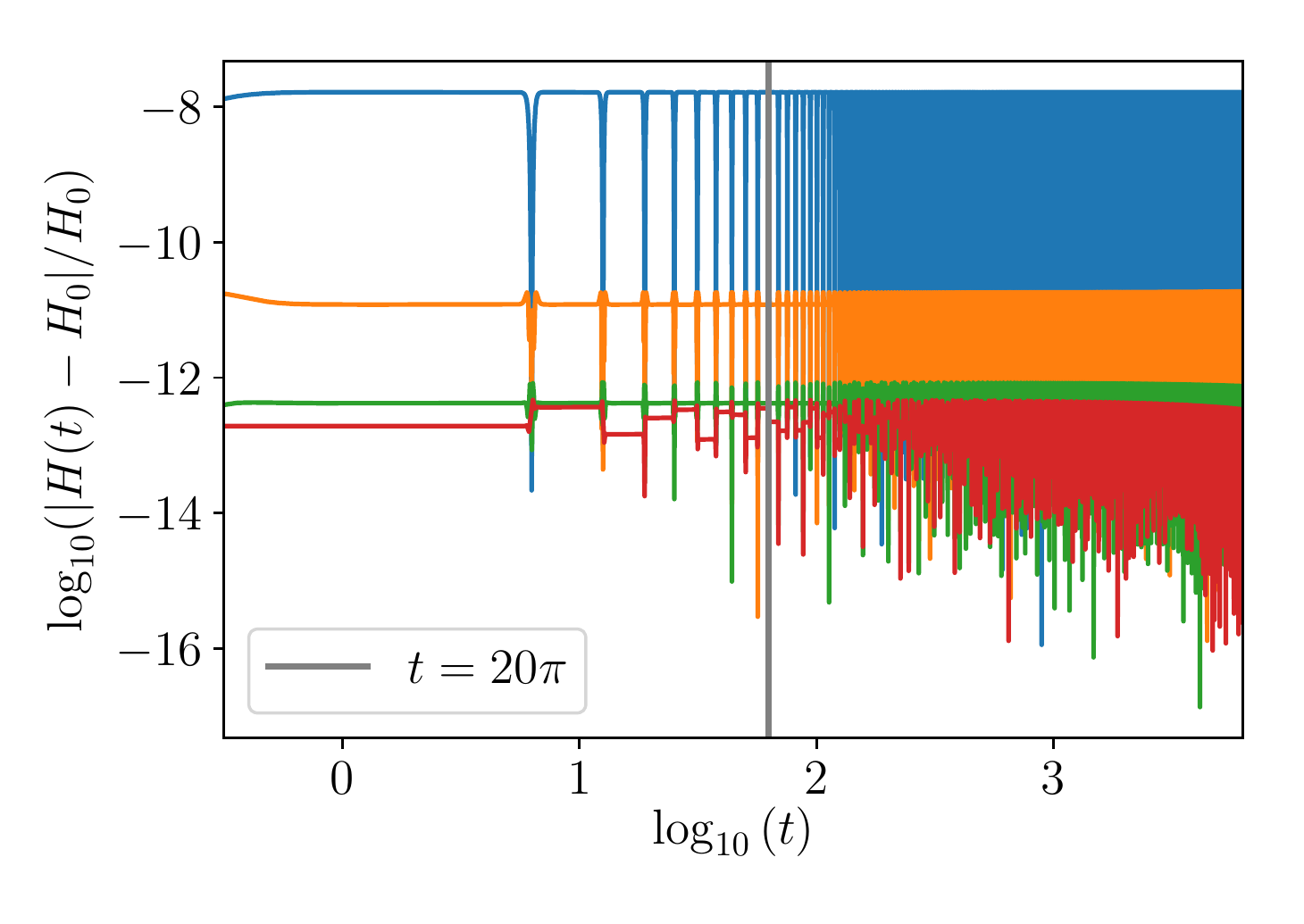}
    \caption{\small Left: Relative error in energy vs. number of evaluations of the basic scheme $\mathcal{S}^{[4]}$ (blue) for the different $T$-methods: 
    $T^{(1)}$ (orange), $T^{(2)}$ (green) and $T^{(3)}$ (red), in the interval $t \in [0, 20 \pi]$. Right: Evolution of this error along the integration when
    $t \in [0, 2000 \pi]$. In this case the step size is chosen so that all schemes involve the same number of evaluations of the basic method.}
    \label{fig:Tkep}
  \end{center}
\end{figure}

We take $e = 0.6$, integrate until the final time $t_f = 20 \pi$ with the basic splitting method $\psis_h^{[4]}$ given by (\ref{s4sim})
and schemes $T_h^{(k)}$, with $k=1,2,3$ for several time steps and then we compute the average error in energy along the integration
interval. Figure \ref{fig:Tkep} (left) shows this error as a function of the number of evaluations of the basic scheme $\psis_h^{[4]}$.
The diagram clearly exhibits the order of convergence of each method: order 4 for $\psis_h^{[4]}$, and orders 6, 8 and 10 for
$T_h^{(1)}$, $T_h^{(2)}$ and $T_h^{(3)}$, respectively.

In the right panel we show the long-time behavior of the error in energy for each method when the step size is chosen so that all of them
involve the same computational cost. We see that the error in energy is almost constant for $t \le 2000 \pi$, as is the case for
symplectic integrators. In other words, the lack of symplecticity at order $h^{12}$ has no effect in this integration interval.
In addition, the scheme $T_h^{(3)}$ provides the smaller error.

\paragraph{A linear parabolic equation.}
Our second example concerns the linear equation in one-dimension
\begin{equation} \label{eq:lrd}
\frac{\partial }{\partial t}u(x,t) =  \frac{\partial^2}{\partial x^2}u(x,t) + V(x) u(x,t),  \qquad u(x,0) =  \sin(2 \pi x),
\end{equation}
with periodic boundary conditions in the space domain $[0,1]$. We take 
$V(x)=8+ 4 \sin(2 \pi x)$ and partition the interval $[0,1]$ into $N$ parts of length $\Delta x = 1/N$, so that the vector
$U = (U_0, \ldots, U_{N-1})^T \in \mathbb{R}^N$ is formed, with $U_j = u(x_j, t)$ and $x_j = j/N$, $j=0,1,\ldots, N-1$.
If a Fourier spectral collocation method is used, 
we end up with the $N$-dimensional linear ODE
\begin{equation} \label{eq:problem0}
\frac{dU}{dt} = A \, U + B \, U,
\end{equation}
where $B = \mathrm{diag} (V(x_0), \ldots, V(x_{N-1}))$ and $A$ is a (full) differentiation matrix related with the second derivative $\partial_{xx}$.
The splitting here corresponds to solving separately the systems $\dot{U} = A \, U$ and
$\dot{U} = B \, U$. Notice that, since $B$ is diagonal, then
\[
  (\e^{h B} U)_j = \e^{h V(x_j)} U_j
\]
and only requires the computation of $N$ multiplications. On the other hand, $A \, U = \mathcal{F}^{-1} D_A \mathcal{F} \, U$, where $\mathcal{F}$ and
$\mathcal{F}^{-1}$ are the forward and backward discrete Fourier transform, and $D_A$ is again diagonal \cite{trefethen00smi}. In consequence,
\[
  \e^{h A} \, U = \mathcal{F}^{-1} \e^{h D_A} \mathcal{F} \, U,
\]
requiring $\mathcal{O}(N \log N)$ operations when the transformation $\mathcal{F}$ (and its inverse) is computed with the fast Fourier transform
(FFT) algorithm.

\begin{figure}[!h]
  \begin{center}
    \includegraphics[width=12cm]{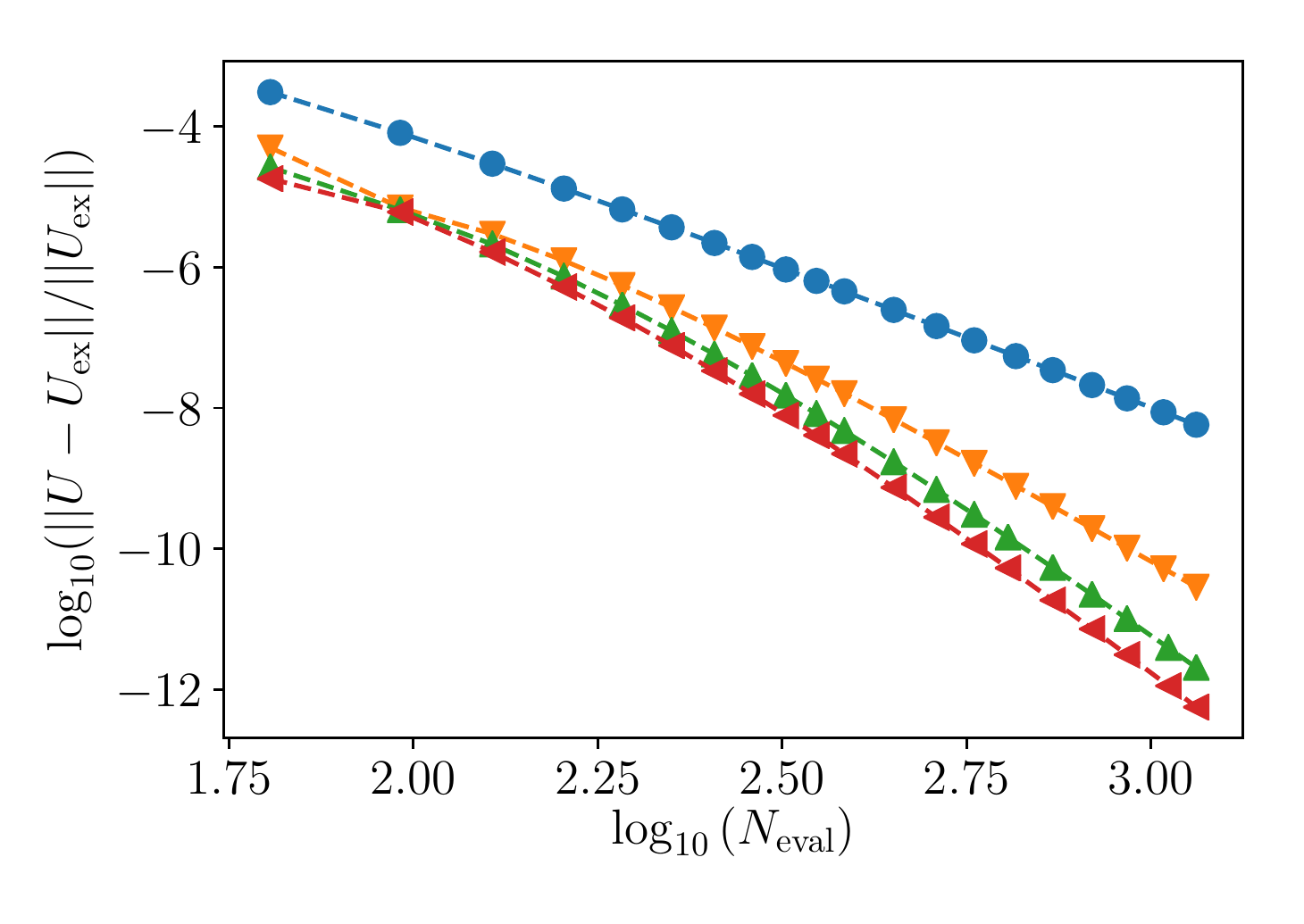}
    \caption{Error vs. number of evaluations of the basic method  (\ref{s4sim}) obtained by schemes $T_h^{(1)}$ (orange), $T_h^{(2)}$ (green) and $T_h^{(3)}$ (red).
    The blue line corresponds to $\psis_h^{[4]}$.}
    \label{fig:Tpar}
  \end{center}
\end{figure}

We take $N=128$ and integrate until $t_f = 1$, where we compute the relative error $\|U-U_{ex}\|/\|U_{ex}\|$ with each method $T_h^{(k)}$, $k=1,2,3$, in addition to the
basic scheme (\ref{s4sim}). The `exact' solution $U_{ex}$ is taken as the output of the 8th-order composition method \textbf{P8S15} of \cite{blanes13oho}. 
The corresponding efficiency diagram is shown in Figure \ref{fig:Tpar}, where the same notation is used for the curves depicted. Here also the higher degree
integrators provide the best efficiency.

\section{Discussion}

\subsection{$T$-methods and $R$-methods}

Methods $T_h^{(k)}$ have indeed close similarities with the compositions $R_h^{(k)}$ (\ref{foit}) previously analyzed in
\cite{casas21cop}: not only their starting point is the same (the basic time-symmetric method $\psis_h^{[2n]}$), but one has in addition
$T_h^{(1)} = R_h^{(1)}$ and also the same coefficients $\gamma^{[2m]}$ defined in (\ref{5c}) enter into their formulation. Finally, they have
the same preservation properties. There is, however, a fundamental difference: whereas $T$-methods are linear combinations of symmetric-conjugate compositions
only, this is not the case of $R$-methods, and in fact schemes $R_h^{(k)}$ involve a much larger number of compositions. This can be clearly
seen by writing explicitly the expression of $R_h^{(2)}$:
\begin{eqnarray} \label{eq.r2}
  R_{h}^{(2)}  & = & \frac{1}{8} \Big( (\gamma^{[2n+2]} \gamma^{[2n]}, \gamma^{[2n+2]} \bar{\gamma}^{[2n]}, \bar{\gamma}^{[2n+2]} \gamma^{[2n]},
         \bar{\gamma}^{[2n+2]} \bar{\gamma}^{[2n]} ) \nonumber \\
   &    & \quad  +   (\gamma^{[2n+2]} \bar{\gamma}^{[2n]}, \gamma^{[2n+2]} \gamma^{[2n]}, \bar{\gamma}^{[2n+2]} \bar{\gamma}^{[2n]},
         \bar{\gamma}^{[2n+2]} \gamma^{[2n]} ) \\
   &   & \quad  +   (\gamma^{[2n+2]} \gamma^{[2n]}, \gamma^{[2n+2]} \bar{\gamma}^{[2n]}, \bar{\gamma}^{[2n+2]} \bar{\gamma}^{[2n]},
         \bar{\gamma}^{[2n+2]} \gamma^{[2n]} )   \nonumber \\  
   &    & \quad +  (\gamma^{[2n+2]} \bar{\gamma}^{[2n]}, \gamma^{[2n+2]} \gamma^{[2n]}, \bar{\gamma}^{[2n+2]} \gamma^{[2n]},
         \bar{\gamma}^{[2n+2]} \bar{\gamma}^{[2n]} ) \nonumber \\
   &  &  \;\; + \mathrm{c.c.}      \Big),   \nonumber
\end{eqnarray}
whereas $R_{h}^{(3)}$ is the sum of 64 compositions  containing 8 basic schemes with weights $\gamma^{[2n+4]} \gamma^{[2n+2]} \gamma^{[2n]}$,
$\bar{\gamma}^{[2n+4]} \gamma^{[2n+2]} \gamma^{[2n]}$, etc. plus their complex conjugate divided by 128. In general,
$R_{h}^{(k)}$ involves the sum of $2^{2^k-2}$ compositions of $2^k$ appropriately weighted basic schemes:
\[
  R_{h}^{(k)} = \frac{1}{2^{2^k -1}} \sum_{j=1}^{2^{2^k-2}} \Big( (\alpha_{j_{2^k}},  \ldots,  \alpha_{j_1})  + \mbox{c.c.} \Big),
\]
where  $\alpha_{j_i}$ are products of the $k$ coefficients   $\gamma^{[2n]}, \ldots, \gamma^{[2(n+k-1)]}$ and their complex
conjugate. This should be compared with the $T$-methods: in general, $T_{h}^{(k)}$ 
involves the sum of $2^{k-1}$ compositions of $2^k$ basic schemes. In either case, the computation of the complex conjugate part can be
avoided just by taking the real part, with no extra evaluations of $\psis_h^{[2n]}$.

\begin{table}[!h]
  \renewcommand\arraystretch{1.2}
  \begin{center}
    \begin{tabular}{c|c|c|c}
 $k$ & $R$ (explicit) & $R$ (recursive) & $T$ (explicit) \\ \hline
 $1$ & $2$                &		$2$	    &  $2$  \\
 $2$ & $16$		&		$8$	    &  $8$ \\
 $3$ & $512$		&		$32$	    &  $32$ \\
  $\vdots$ & $\vdots$&$\vdots$&$\vdots$ \\
 $m$ & $2^m \cdot 2^{2^m-2}$ & 	$2^m \cdot 2^{m-1}$	    & $2^m \cdot 2^{m-1}$ 
\end{tabular} 
\caption{\small Number of basic maps $\psis_h^{[2n]}$ necessary to compute when formulating $R$- and $T$-methods explicitly or recursively (in the case of
$R$-methods).}
    \label{tau:costR}
  \end{center}
\end{table}

These numbers are collected in Table \ref{tau:costR}, when schemes $R_h^{(k)}$  (second column) and $T_h^{(k)}$ (last column) are
formulated explicitly. Of course, a recursive implementation of $R$-methods by applying the procedure (\ref{foit}) turns out to be more
efficient. In that case the required computational effort, measured as the number of basic schemes, is shown in the third column of the table. Again,
in this case we only have to compute the real part in the last iteration.


In view of the number of basic maps required by the recursive implementation of $R$-methods and the explicit formulation (\ref{T.methods}) 
of $T$-methods, it is natural to ask what are the advantages (if any) of the later schemes with respect to the former ones. In this respect,
one should take into account that both explicit formulations (\ref{T.methods}) and (\ref{eq.r2}) are directly amenable to parallelization, whereas
this is less obvious for the recursion (\ref{foit}). 

If one has a computer with, say, $2\ell$ threads, it is easy to estimate the effective number of
evaluations of $\psis_h^{[2n]}$ both for $R$- and $T$-methods. Thus, for   $R_{h}^{(k)}$ one has:
\begin{itemize}
 \item if $\ell \le 2^k -2$ then the number of evaluations is $2^k$;
 \item if $\ell > 2^k -2$ then the number of evaluations is $2^k \cdot 2^{2^k-2-\ell}$,
\end{itemize}
whereas this number is considerably reduced for schemes $T_{h}^{(k)}$:
\begin{itemize}
 \item if $\ell \le k-1$ then the number of evaluations is $2^k$;
 \item if $\ell > k-1$ then the number of evaluations is $2^k \cdot 2^{k-1-\ell}$.
\end{itemize}
In Table \ref{cost.paralell} we collect these numbers for the first values of $k$ in the particular case of $2^{2} = 4$ and 
$2^{5} = 32$ threads. We see that the implementation of the explicit expression
of the $R$-methods is more advantageous than the recursive procedure already with a relatively small number of threads, and that, in any
case, $T$-methods require less computational effort.
  
\begin{table}[!h]
  \renewcommand\arraystretch{1.2}
  \begin{center}
    \begin{tabular}{|c|c|c|c|c|}
     & \multicolumn{2}{c|}{ $4$ threads} & \multicolumn{2}{c|}{ $32$ threads} \\ \hline
  $k$  & $R$  &  $T$   & $R$  &  $T$\\ \hline
 $1$ & $2$   & $2$    & $2$  & $2$ \\
 $2$ & $4$   & $4$   & $4$  & $4$ \\
 $3$ & $128$ &	 $8$  & $16$ & $8$ \\ \hline
\end{tabular} 
\caption{\small Effective number of evaluations of the basic map $\psis_h^{[2n]}$ when the corresponding $R$- and $T$-method is implemented in
parallel with 4 and 32 threads.}
    \label{cost.paralell}
  \end{center}
\end{table}

\begin{figure}[!h]
  \begin{center}
    \includegraphics[scale=1.00]{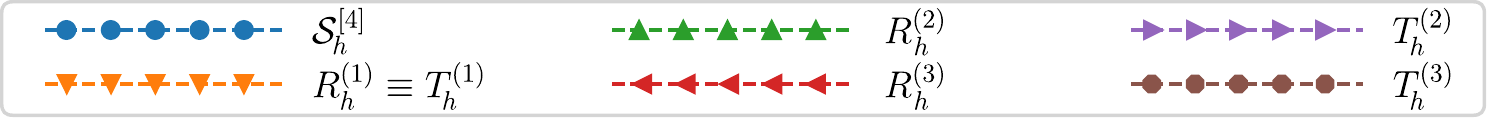}\\
    \includegraphics[scale=0.55]{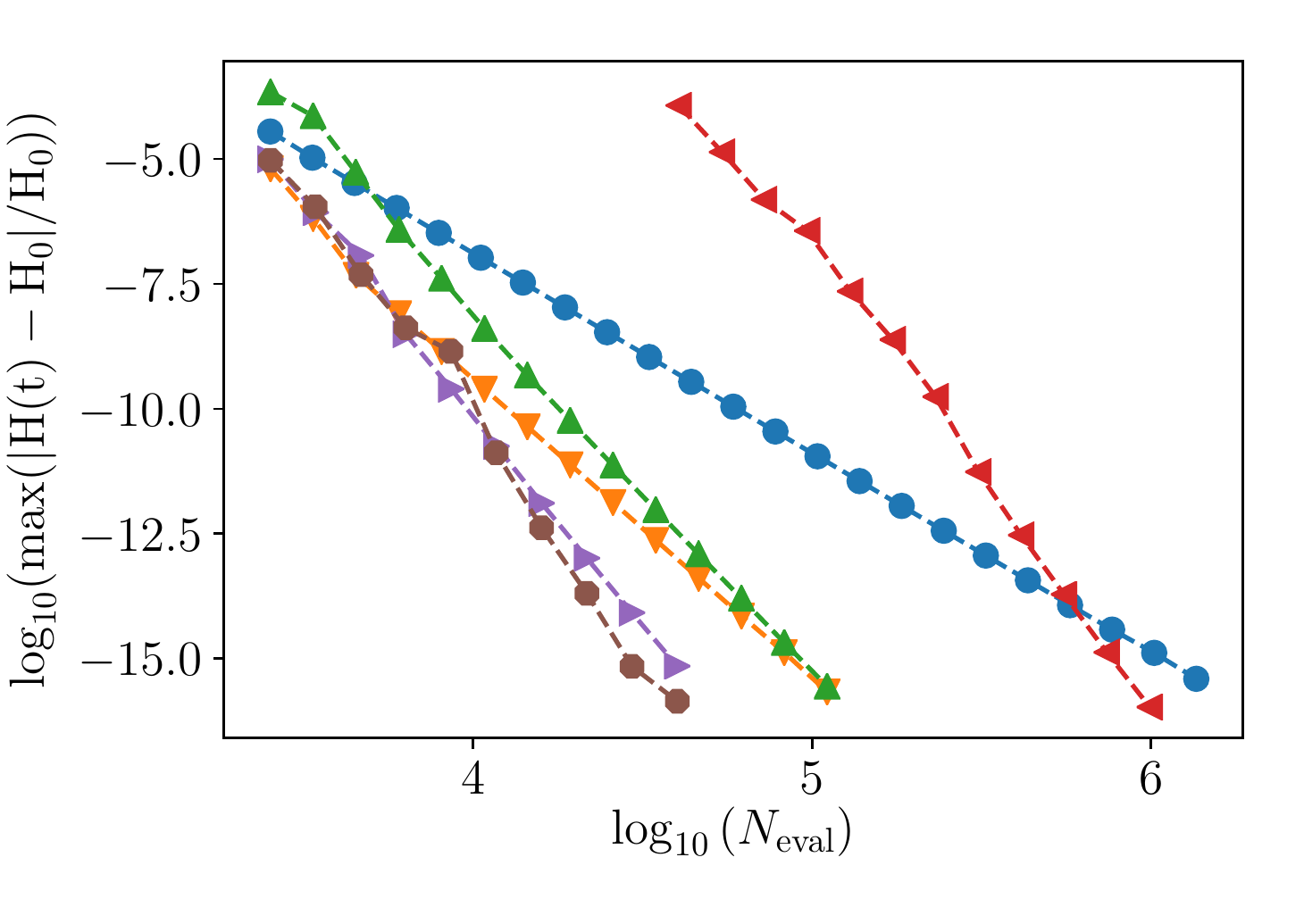}
    \includegraphics[scale=0.55]{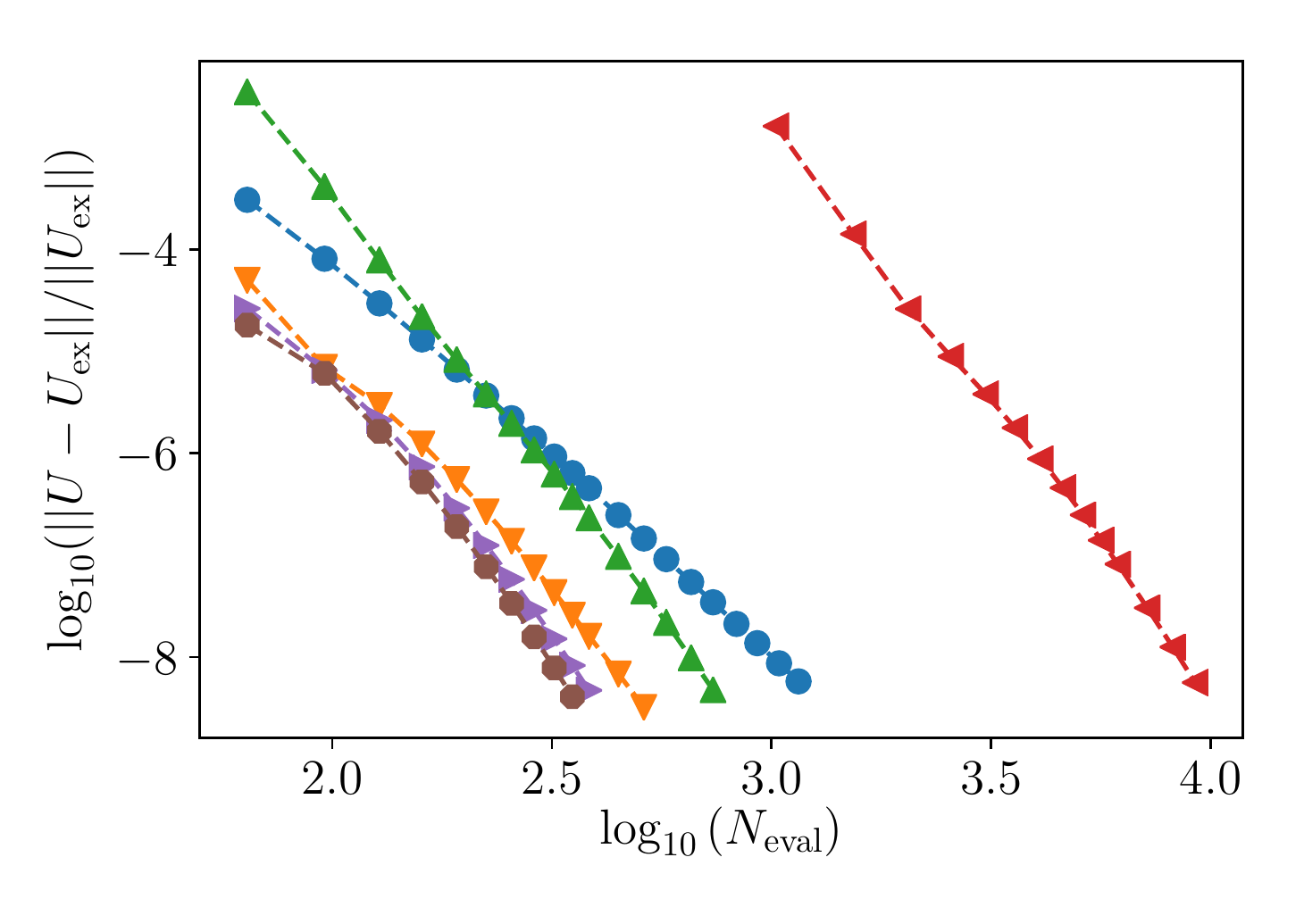}
    \caption{\small Efficiency diagram for $R$- and $T$-methods in a computer accepting 4 threads without loss of performance. Left: Kepler problem. Right: linear
    parabolic equation.}
    \label{fig:4threads}
  \end{center}
\end{figure}

To better illustrate this issue, we next compare the efficiency of the different methods when implemented on a computer able to
execute 4 threads without loss of performance. The corresponding results are displayed in Figure \ref{fig:4threads} for the Kepler problem (left) and
the linear parabolic  equation (\ref{eq:lrd}) (right). The gain in efficiency of the new schemes is clearly visible.

\begin{figure}[!h]
  \begin{center}
    \includegraphics[scale=1.00]{legend3.pdf}\\
    \includegraphics[scale=0.55]{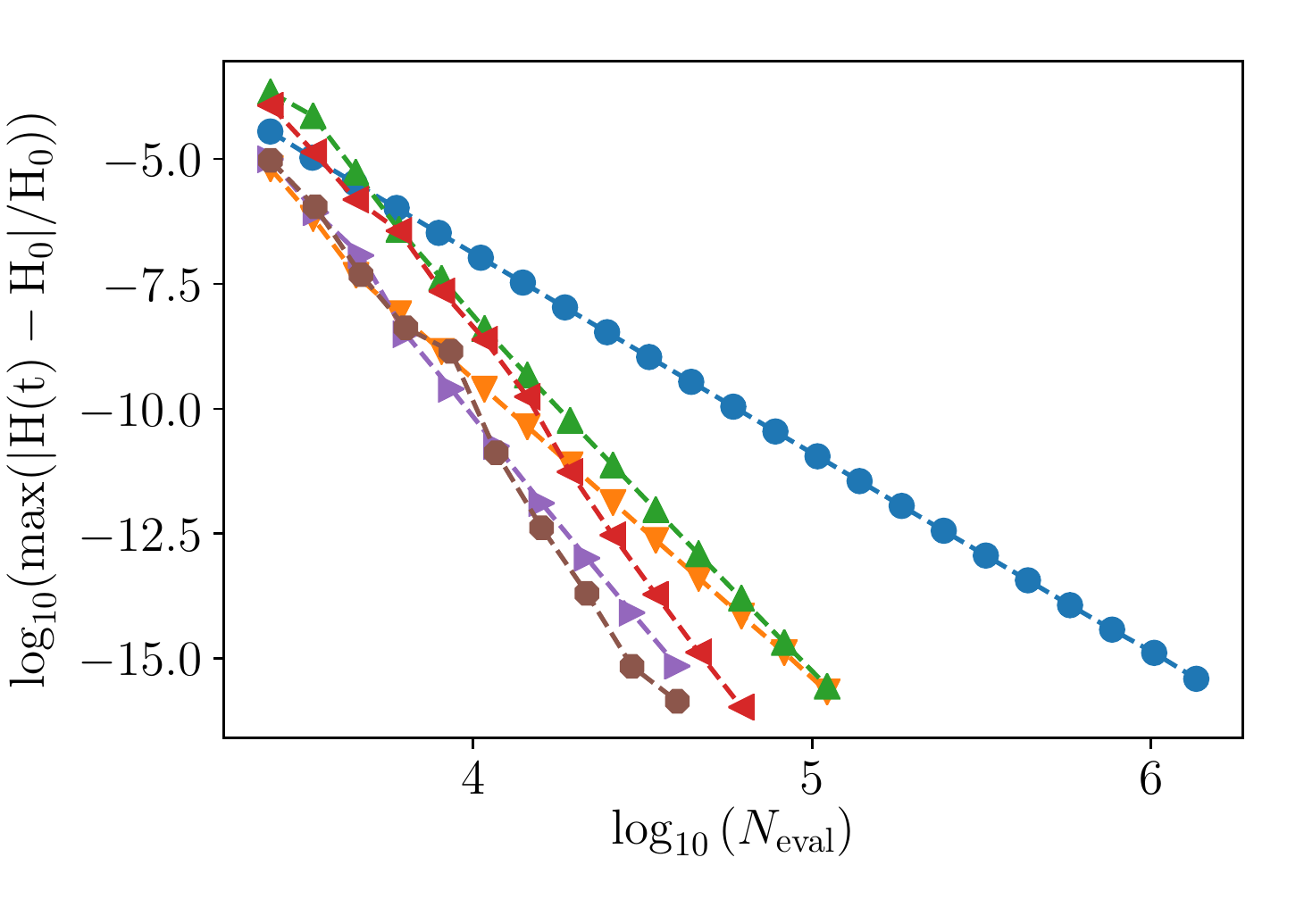}
     \includegraphics[scale=0.55]{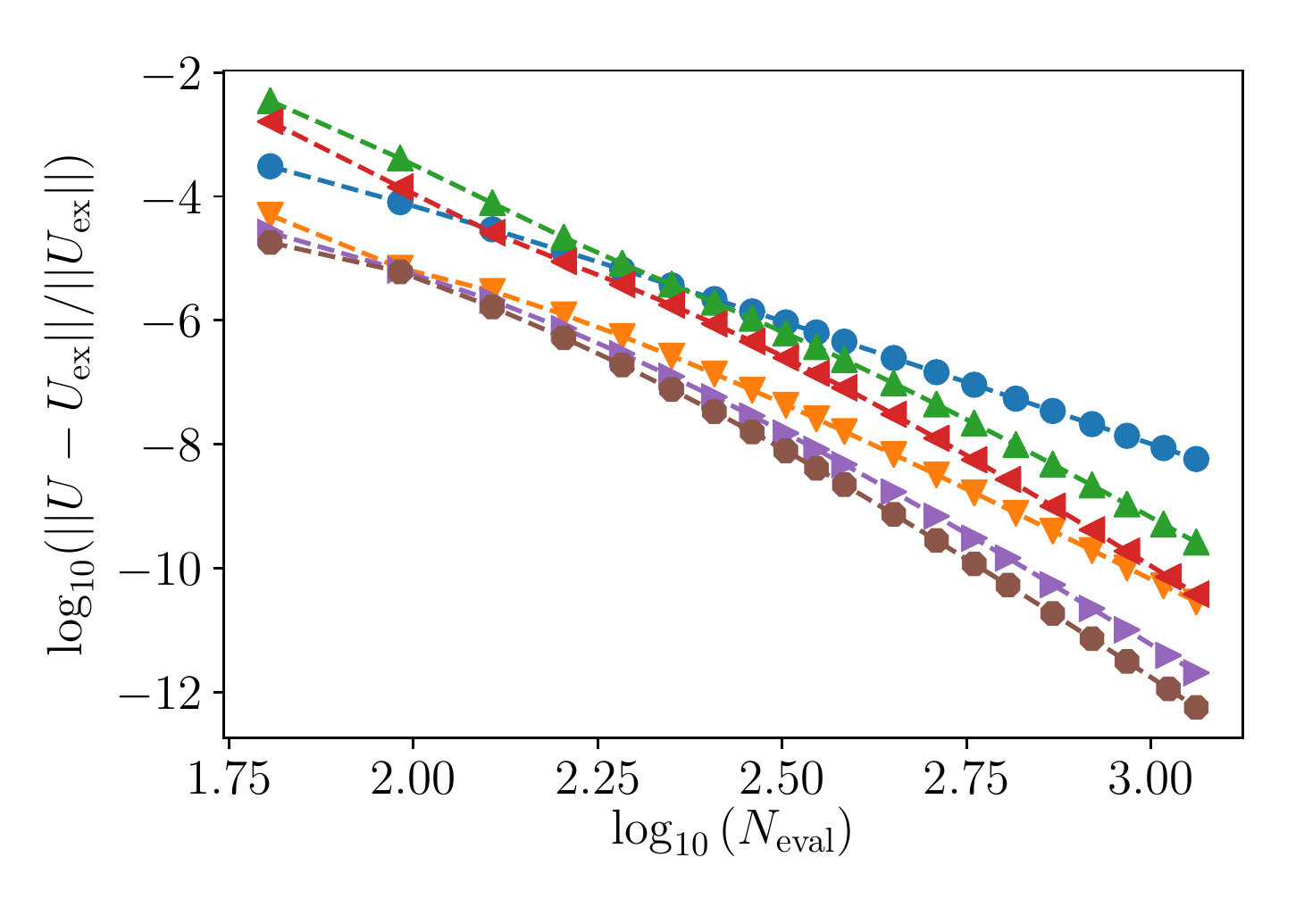}
    \caption{\small Efficiency diagram for $R$- and $T$-methods when the effective number of evaluations of the basic scheme is the same ($2^k$ in
    both cases). Left: Kepler problem. Right: linear parabolic equation.}
    \label{fig:32threads}
  \end{center}
\end{figure}

Even in the case when one could run the schemes on a machine such that the effective number of evaluations of both  $R_{h}^{(k)}$ and $T_{h}^{(k)}$
is the same, i.e., $2^k$ in both cases,  the latter turn out to be more efficient. This is clearly visible in Figure \ref{fig:32threads}, obtained again by applying the previous schemes to the
Kepler problem (left) and the linear parabolic equation (right). 

\begin{figure}[!h]
  \begin{center}
    \includegraphics[scale=1.00]{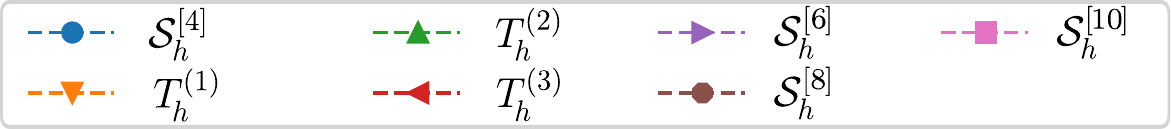}\\
    \includegraphics[width=12cm]{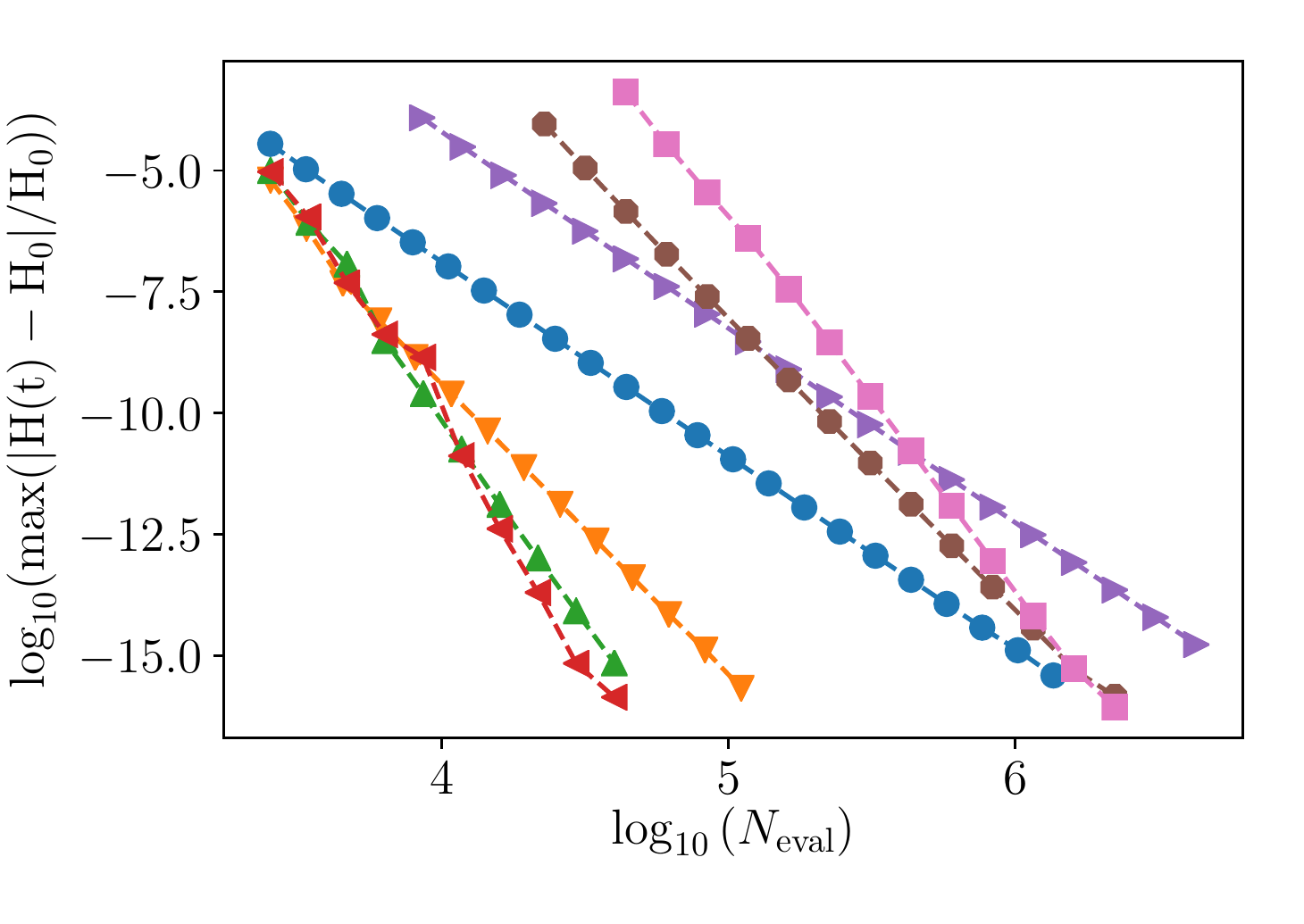}
    \caption{\small Error in energy vs. number of basic maps $\psis_h^{[4]}$ for $T$-methods in comparison with schemes obtained by triple-jump for the
    Kepler problem.}
    \label{fig:TJ}
  \end{center}
\end{figure}

Finally, it is also illustrative to compare the efficiency of the new $T$-methods with the standard triple-jump procedure, Eqs. (\ref{eq:tj})-(\ref{real.1}), both applied
to the same basic scheme (\ref{s4sim}). Thus, in Figure \ref{fig:TJ} we depict the results achieved by projecting $\psis_h^{[6]}$, $\psis_h^{[8]}$, and 
$\psis_h^{[10]}$ at each step,
together with $T_{h}^{(k)}$, $k=1,2,3$ for the Kepler problem with the same parameters and final time $t_f = 20 \pi$. Here the effective number of evaluations of the
basic scheme has been taken as $2^k$ for $T$-methods and $3^k$ for triple-jump.
Not surprisingly, the new schemes turn out to be much more efficient.

\subsection{Concluding remarks}

The standard triple-jump procedure is a popular technique that allows one to construct numerical integrators for differential equations of arbitrarily high order
by composition of a basic integrator of low order. It has nevertheless certain limitations: the number of basic maps grows rapidly
with the order, and the main error terms are quite large in comparison with other specially built integrators. Moreover, they involve some negative coefficients
when the order $r \ge 3$, so that the resulting schemes cannot be used in particular when the initial value problem (\ref{ivp}) results from the space
discretization of a parabolic partial differential equation involving the Laplace operator.  In this context it is quite natural to explore whether it is still possible using
the triple-jump technique (\ref{eq:tj}), but with the complex coefficients furnished by (\ref{cco1}) as long as their real part is positive. It has been established that
this is indeed the case, although once again they require an exceedingly large number of basic methods. For this reason, other
alternatives for constructing high-order composition methods have also been proposed \cite{blanes13oho,castella09smw,hansen09hos}.  
Among them, the class of schemes (\ref{R.1}) possess
some special features: starting from a time-symmetric basic scheme $\psis_h^{[2n]}$ of order $2n$, it is possible to construct recursively methods of order
$2n+2k$, $k=1,2,\ldots$ that are still time-symmetric up to order $4n+3$. Moreover, if the differential equation in (\ref{ivp}) has some qualitative properties (such
as symplecticity or volume preservation) then these properties are still shared by the numerical solution up to order $4n+3$ \cite{casas21cop}. 

Methods (\ref{foit}) are based on the simple symmetric-conjugate composition (\ref{sc.1}). As shown in \cite{blanes21osc}, 
sym\-me\-tric-conjugate composition methods still
possess remarkable preservation properties when projected on the real axis at each integration step, and so it makes sense to consider more general linear
combinations of methods within this class. The corresponding analysis has been carried out here, and as a result we have built a new class of schemes that
essentially have the same preservation properties as methods (\ref{foit}), but requiring a much reduced computational cost. In addition, these methods are
particularly well suited for their parallel implementation. The examples included show a significant improvement in efficiency with respect to schemes
(\ref{foit}) and those obtained by applying the triple-jump procedure.

\subsection*{Acknowledgements}
This work has been funded by  
Ministerio de Ciencia e Innovaci\'on (Spain) through project PID2019-104927GB-C21 (AEI/FEDER, UE) and by Universitat Jaume I (UJI-B2019-17). 
A.E.-T. has been additionally supported by the predoctoral contract BES-2017-079697 (Spain). 


\begin{thebibliography}{10}

\bibitem{arnold89mmo}
{\sc V.~Arnold}, {\em Mathematical {M}ethods of {C}lassical {M}echanics},
  Springer-Verlag, {S}econd~ed., 1989.

\bibitem{aubry98psr}
  {\sc A.~Aubry and P.~Chartier}, {\em Pseudo-symplectic Runge--Kutta methods},
  BIT Num. Math., 38 (1998), pp.~439--461.
  
\bibitem{bandrauk91ies}
{\sc A.~Bandrauk and H.~Shen}, {\em Improved exponential split operator method
  for solving the time-dependent {S}chr\"{o}dinger equation}, Chem. Phys.
  Lett., 176 (1991), pp.~428--432.

\bibitem{blanes99eos}  
{\sc S.~Blanes, F.~Casas, and J. Ros}, {\em Extrapolation of symplectic integrators},
Celest. Mech. \& Dyn. Astr., 75 (1999), pp.~149--161.

\bibitem{blanes05otn}
{\sc S.~Blanes and F.~Casas}, {\em On the necessity of negative coefficients
  for operator splitting schemes of order higher than two}, Appl. Numer. Math.,
  54 (2005), pp.~23--37.

\bibitem{blanes16aci}
{\sc S.~Blanes and F.~Casas}, {\em A {C}oncise
  {I}ntroduction to {G}eometric {N}umerical {I}ntegration}, {CRC} Press, 2016.

\bibitem{blanes21osc}
{\sc S.~Blanes, F.~Casas, P.~Chartier, and A.~Escorihuela-Tom\`as}, {\em On
  symmetric-conjugate composition methods in the numerical integration of
  differential equations}, Tech. Rep. 2101.04100, arXiv, 2021.

\bibitem{blanes13oho}
{\sc S.~Blanes, F.~Casas, P.~Chartier, and A.~Murua}, {\em Optimized high-order
  splitting methods for some classes of parabolic equations}, Math. Comput., 82
  (2013), pp.~1559--1576.

\bibitem{blanes08sac}
{\sc S.~Blanes, F.~Casas, and A.~Murua}, {\em Splitting and composition methods
  in the numerical integration of differential equations}, Bol. Soc. Esp. Mat.
  Apl., 45 (2008), pp.~89--145.

\bibitem{blanes10smw}
{\sc S.~Blanes, F.~Casas, and A.~Murua}, {\em Splitting methods
  with complex coefficients}, Bol. Soc. Esp. Mat. Apl., 50 (2010), pp.~47--61.

\bibitem{casas21cop}
{\sc F.~Casas, P.~Chartier, A.~Escorihuela-Tom\`as, and Y.~Zhang}, {\em
  Compositions of pseudo-symmetric integrators with complex coefficients for
  the numerical integration of differential equations}, J. Comput. Appl. Math.,
  381 (2021), p.~113006.

\bibitem{castella09smw}
{\sc F.~Castella, P.~Chartier, S.~Descombes, and G.~Vilmart}, {\em Splitting
  methods with complex times for parabolic equations}, BIT Numer. Math., 49
  (2009), pp.~487--508.

\bibitem{chambers03siw}
{\sc J.~Chambers}, {\em Symplectic integrators with complex time steps},
  Astron. J., 126 (2003), pp.~1119--1126.

\bibitem{chan00eos}
{\sc R. Chan and A. Murua}, {\em Extrapolation of symplectic methods for Hamiltonian  
problems}, Appl. Numer. Math., 34 (2000), pp.~189--205.


\bibitem{goldman96nno}
{\sc D.~Goldman and T.~Kaper}, {\em $n$th-order operator splitting schemes and
  nonreversible systems}, SIAM J. Numer. Anal., 33 (1996), pp.~349--367.

\bibitem{hairer06gni}
{\sc E.~Hairer, C.~Lubich, and G.~Wanner}, {\em Geometric {N}umerical
  {I}ntegration. {S}tructure-{P}reserving {A}lgorithms for {O}rdinary
  {D}ifferential {E}quations}, Springer-Verlag, {S}econd~ed., 2006.

\bibitem{hansen09hos}
{\sc E.~Hansen and A.~Ostermann}, {\em High order splitting methods for
  analytic semigroups exist}, BIT Numer. Math., 49 (2009), pp.~527--542.

\bibitem{mclachlan02sm}
{\sc R.~McLachlan and R.~Quispel}, {\em Splitting methods}, Acta Numerica, 11
  (2002), pp.~341--434.

\bibitem{sanz-serna94nhp}
{\sc J.~Sanz-Serna and M.~Calvo}, {\em Numerical {H}amiltonian {P}roblems},
  Chapman {\&} Hall, 1994.

\bibitem{sheng89slp}
{\sc Q.~Sheng}, {\em Solving linear partial differential equations by
  exponential splitting}, IMA J. Numer. Anal., 9 (1989), pp.~199--212.

\bibitem{suzuki90fdo}
{\sc M.~Suzuki}, {\em Fractal decomposition of exponential operators with
  applications to many-body theories and {M}onte {C}arlo simulations}, Phys.
  Lett. A, 146 (1990), pp.~319--323.

\bibitem{suzuki91gto}
{\sc M.~Suzuki}, {\em General theory of
  fractal path integrals with applications to many-body theories and
  statistical physics}, J. Math. Phys., 32 (1991), pp.~400--407.

\bibitem{trefethen00smi}
{\sc L.N.~Trefethen}, {\em Spectral Methods in MATLAB},
  SIAM, 2000.

\bibitem{yoshida90coh}
{\sc H.~Yoshida}, {\em Construction of higher order symplectic integrators},
  Phys. Lett. A, 150 (1990), pp.~262--268.

\end{thebibliography}

\end{document}